\documentstyle[twoside,11pt]{article}
  \newcommand{\bref}[5]{\noindent\parbox[t]{0.8cm}{#1}\parbox[t]{15.2cm}{{#2}{\it #3}{\bf #4}#5}\par}
   
  \def\vp{\varepsilon}
  
\begin{document}
  \title{\bf {Some criteria  for integer sequences pair  being realizable by a graph}
  \thanks{This work is supported by the National Natural Science Foundation of
  China(NSFC11921001,11601108) and the National Key Research and Development Program of
  China(2018YFA0704701).  }}
\author{ $ {\mbox {Jiyun \,Guo}}$, $ {\mbox {Miao \,Fu*}}$, ${\mbox {Jun \,Wang}}$ \\
{\small  College of Mathematics, Tianjin University,}
{\small Tianjin, 300072, P. R. China } \\
 {\small Email: fumiao1119@163.com} }
\date{}
\maketitle
\begin{center}
\parbox{0.9\hsize}
{\small {\bf Abstract.}\ \ Let $A=(a_1,\ldots,a_n)$ and
$B=(b_1,\ldots,b_n)$ be two sequences of
 nonnegative integers with $a_i \le
b_i$ for $1\le i\le n$. The pair $(A;B)$ is said to be realizable by
a graph if there exists a simple graph $G$ with vertices
$v_1,\ldots, v_n$ such that $a_i\le d_G(v_i)\le b_i$ for $1\le i\le
n$. Let $\preceq$ denote the lexicographic ordering on $Z\times Z:$
$(a_{i+1},b_{i+1})\preceq (a_i,b_i)\Longleftrightarrow
[(a_{i+1}<a_i)\vee ((a_{i+1}=a_i)\&(b_{i+1}\le b_i))]$. We say that
the sequences $A$ and $B$ are in good order if
$(a_{i+1},b_{i+1})\preceq (a_i,b_i)$. In this paper, we consider the
generalizations of six classical characterizations on sequences pair
due to Berge, Ryser et al. and present related results.

{\bf Keywords.}\ \ Criteria; Realizable; Characterization; Good
order.

{\bf Mathematics Subject Classification(2000):} 05C07.}
\end{center}

\section*{1. Introduction}
\hskip\parindent Throughout this paper, we consider only undirected
graphs without loops or parallel edges.

A non-increasing sequence  $\pi=(d_1,\ldots, d_n)$ of nonnegative
integers with even sum is said to be realizable if there exists a
graph $G$ having $\pi$ as one of its degree sequences, in which case
we call $G$ a realization of $\pi$. In [13], seven criteria for such
an integer sequence being realizable are listed, let us start with
Berge criterion. Before presenting it, we first introduce Berge
sequence.

Define Berge sequence $\bar{\pi}=(\bar{d}_1,\ldots,\bar{d}_n)$ as
follows: For $i\in \{1,\ldots,n\}$, $\bar{d}_i$ is the $i$th column
sum of the $(0,1)$-matrix, which the $d_k$(has for each $k$) leading
terms in row $k$ equal to 1 except for the $(k,k)th$ term that is 0
and also the remaining entries are 0's. For example, if $d_1=4$,
$d_2=2$, $d_3=2$, $d_4=2$, $d_5=1$, and the $(0,1)$-matrix becomes
$$\left({\begin{array}{ccccc}
{{0}}&{{1}}&{{1}}&{{1}}&{{1}}\\
{{1}}&{{0}}&{{1}}&{{0}}&{{0}}\\
 {{1}}&{{1}}&{{0}}&{{0}}&{{0}}\\
{{1}}&{{1}}&{{0}}&{{0}}&{{0}}\\
{{1}}&{{0}}&{{0}}&{{0}}&{{0}}\\
\end{array}}\right).$$

{\bf Theorem 1.1.}(Berge [1]). $\pi=(d_1,\ldots, d_n)$ is realizable
if and only if $$\sum\limits_{i=1}^k d_i\le \sum\limits_{i=1}^k
\bar{d}_i\mbox{ for each $k$ with $1\le k\le n$.}$$

A sequence $(h_1,\ldots,h_m;l_1,\ldots,l_n)$ is called bipartite
realizable iff there is a bipartite graph such that one component
has degree sequence $(h_1,\ldots,h_m)$ and the other one has
$(l_1,\ldots,l_n).$ Define $g=max\{i|d_i\ge i\}$ and
$\tilde{d}_i=d_i+1$ if $i\in <g>(=\{1,\ldots,g\})$ and
$\tilde{d}_i=d_i$ otherwise. The following is Ryser criterion.

{\bf Theorem 1.2.}(Ryser [2]). $\pi=(d_1,\ldots, d_n)$ is realizable
if and only if
$(\tilde{d_1},\ldots,\tilde{d_n};\tilde{d_1},\ldots,\tilde{d_n})$ is
bipartite realizable.

{\bf Theorem 1.3.}(Erd$\ddot{ o}$s and Gallai [5]).
$\pi=(d_1,\ldots, d_n)$ is realizable if and only if

$$\sum_{i=1}^k d_i\le  k(k-1)+\sum_{i=k+1}^n
min\{d_i,k\}\mbox{ for each $k$ with $1\le k\le n$.}$$

{\bf Theorem 1.4.}(Fulkerson et al. [6]).  $\pi=(d_1,\ldots, d_n)$
is realizable if and only if

$$\sum_{i=1}^k d_i\le  k(n-m-1)+\sum_{i=n-m+1}^n
d_i \mbox{ for each $k$ with $1\le k\le n$,}$$ where $m\ge 0$ and
$k+m\le n$.

{\bf Theorem 1.5.}(Bollob$\acute{ a}$s [2]).  $\pi=(d_1,\ldots,
d_n)$ is realizable if and only if

$$\sum_{i=1}^k d_i\le  \sum_{i=k+1}^n d_i+
\sum_{i=1}^k min\{d_i,k-1\}\mbox{ for each $k$ with $1\le k\le
n$.}$$

{\bf Theorem 1.6.}(Gr$\ddot{ u}$nbaum [8]).  $\pi=(d_1,\ldots, d_n)$
is realizable if and only if

$$\sum_{i=1}^k max\{k-1,d_i\}\le  k(k-1)+\sum_{i=k+1}^n
d_i \mbox{ for each $k$ with $1\le k\le n$.}$$

H$\ddot{a}$sselbarth [10] also provided such a criterion. They
defined a new sequence $(d^* _1,\ldots,d^* _n)$ as follows: For
$1\le i\le n$, $d^* _i$ is the $i$th column sum of the
$(0,1)$-matrix in which the $d_i$ leading terms in row $i$ are 1's
and the remaining entries are 0's. The sequence $\pi=(d^*
_1,\ldots,d^* _n)$ is called the conjugate of $\pi=(d_1,\ldots,
d_n)$.

{\bf Theorem 1.7.}(H$\ddot{a}$sselbarth [10]).  $\pi=(d_1,\ldots,
d_n)$ is realizable if and only if

$$\sum_{i=1}^k d_i\le \sum_{i=1}^k (d^*_i-1) \mbox{ for each $k$ with $1\le k\le f$,}$$
where $f=max\{i|d_i\ge i\}$.

Motivated by Theorem 1.2, Niessen[12] asked for the extension of
result of Erd$\ddot{o}$s-Gallai to characterizing lists of intervals
$[a_1,b_1],\ldots,[a_n,b_n]$ such that there exists a graph with
vertices  $v_1,\ldots, v_n$ satisfying $a_i\le d(v_i)\le b_i$ for
$1\le i\le n$. Cai et al.[3] presented such a characterization for
$(A;B)$ to be realizable, where $A$ and $B$ are in good
order(defined below). This is a variation of the classical
Lov$\acute{a}$sz's (g,f)-factor theorem[11] in degree sequences,
solves the preceding research problem and generalizes Theorem
1.2(which corresponds to $a_i=b_i=d_i$ for each $i$).

Let $A=(a_1,\ldots,a_n)$ and $B=(b_1,\ldots,b_n)$ be two sequences
of nonnegative integers with $a_i \le b_i$ for $1\le i\le n$. The
pair $(A;B)$ is said to be realizable if there exists a  graph $G$
with vertices $v_1,\ldots, v_n$ such that $a_i\le d_G(x_i)\le b_i$
for $1\le i\le n$. Let $\preceq$ denote the lexicographic ordering
on $Z\times Z:$ $(a_{i+1},b_{i+1})\preceq
(a_i,b_i)\Longleftrightarrow [(a_{i+1}<a_i)\vee
((a_{i+1}=a_i)\&(b_{i+1}\le b_i))]$. We say that the sequences $A$
and $B$ are in good order if $(a_{i+1},b_{i+1})\preceq (a_i,b_i)$.
Cai et al. defined for $t=0,1,\ldots,n$,
$$I_t=\mbox{ $\{i|i\ge t+1$ and $b_i\ge t+1$\}}$$
and
$$\vp(t)=\left\{ \begin{array}{ll}
1, & \mbox{ if $a_i=b_i$ for all $i\in I_t$ and $\sum\limits_{i\in I_t}b_i+t|I_t|$ is odd},\\
0, & \mbox{ otherwise.}\end{array}\right.$$

{\bf Theorem 1.8. }(Cai et al. [3]). \ \ Let $A=(a_1,\ldots,a_n)$
and $B=(b_1,\ldots,b_n)$ be two sequences of nonnegative integers in
good order such that $a_i \le b_i$ for $1\le i\le n$. Then $(A;B)$
is realizable if and only if
$$\sum_{i=1}^t a_i\le t(t-1)+\sum_{i=t+1}^n min\{t,b_i\}-\vp(t)\mbox{ for
each $t$ with $0\le t\le n$.}\eqno{(CDZ)}$$

For the result, Tripathi et al.[7] gave a constructive proof of it.
 Guo and Yin [9] considered the strong situation, Cai and Kang [4]
improved the strong situation and presented a good characterization.

\section*{2. Main Results}
\hskip\parindent  Analogous to the case of Theorem 1.8, we consider
the generalizations of the remaining six classical criteria and give
the relevant conclusions as follows.

{\bf Theorem 2.1.} Let $A=(a_1,\ldots,a_n)$ and $B=(b_1,\ldots,b_n)$
be two sequences of nonnegative integers in good order such that
$a_i \le b_i$ for $1\le i\le n$ and $(\bar{b}_1,\ldots,\bar{b}_n)$
be the Berge sequence of $B$. If the pair $(A;B)$ is realizable by a
graph, then
$$\sum_{i=1}^t a_i\le \sum_{i=1}^t \bar{b}_i\mbox{ for each $t$ with
$0\le t\le n$}.\eqno{(1)}$$

Note that the converse of Theorem 2.1 is not true. For example,
taking $A=(5,4,3,3,3,1)$ and $B=(5,5,3,3,3,1)$. Clearly, $A$ and $B$ satisfy the good
order,
$\bar{B}=(\bar{5},\bar{5},\bar{3},\bar{3},\bar{3},\bar{1})=(5,4,4,3,2,2)$.
 It can be checked that (1) holds for
$t=0,1,2,3,4,5$. However, $(A;B)=(5,4,3,3,3,1;5,5,3,3,3,1)$ is not realizable because $(CDZ)$ does not hold for $t=2$.

Theorem 2.1 and the following theorem  provide a necessary condition
and a sufficient condition, respectively, for the pair $(A;B)$ to be
realizable by a graph with bonds differing by at most one.

{\bf Theorem 2.2.} Let $A=(a_1,\ldots,a_n)$ and $B=(b_1,\ldots,b_n)$
be two sequences of nonnegative integers in good order such that
$a_i \le b_i$ for $1\le i\le n$ and $(\bar{b}_1,\ldots,\bar{b}_n)$
be the Berge sequence of $B$. If the following inequality holds:
$$\sum_{i=1}^t a_i\le \sum_{i=1}^t \bar{b}_i-\varepsilon(t)\mbox{ for each $t$ with
$0\le t\le n$},\eqno{(1')}$$ then the pair $(A;B)$ is realizable by
a graph, where $\varepsilon(t)$ is defined as in Theorem 1.8.

Let $S:=[\tilde{a}_1,\tilde{b}_1],\ldots,[\tilde{a}_n,\tilde{b}_n]$
be a sequence of intervals, where
$\tilde{A}=(\tilde{a}_1,\ldots,\tilde{a}_n)$ and
$\tilde{B}=(\tilde{b}_1,\ldots,\tilde{b}_n)$ are defined as before.
We say that the pair $(S;S)$ is interval realizable if there exists a bipartite graph $F$ with bipartition $\{x_1,\ldots,x_n\}$,
$\{y_1,\ldots,y_n\}$ such that $\tilde{a}_i\le d_F(x_i)\le
\tilde{b}_i$ and $\tilde{a}_i\le d_F(y_i)\le \tilde{b}_i$ for $1\le
i\le n$. The conclusion is as follows:

{\bf Theorem 2.3.} Let $A=(a_1,\ldots,a_n)$ and $B=(b_1,\ldots,b_n)$
be two sequences of nonnegative integers in good order such that
$a_i \le b_i$ for $1\le i\le n$. If the pair $(A;B)$ is realizable,
then $(S;S)$ is interval realizable, where
$S:=[\tilde{a}_1,\tilde{b}_1],\ldots,[\tilde{a}_n,\tilde{b}_n]$.

Note that the converse of Theorem 2.3 is not true. For example,
taking $A=(5,4,3,3,3,1)$ and $B=(5,5,3,3,3,1)$. Clearly, $A$ and $B$ satisfy the
good order, $a_i\le b_i$ for $i\in <6>$.
$\tilde{A}=(\tilde{a}_1,\tilde{a}_2,\tilde{a}_3,\tilde{a}_4,\tilde{a}_5,\tilde{a}_6)=(6,5,4,3,3,1)$,
$\tilde{B}=(\tilde{b}_1,\tilde{b}_2,\tilde{b}_3,\tilde{b}_4,\tilde{b}_5,\tilde{b}_6)=(6,6,4,3,3,1)$ and
$S=6,[5,6],4,3,3,1$. The pair $(S;S)$ is interval realizable because there exists a graph $F$ with bipartition $\{x_1,\ldots,x_n\}$,$\{y_1,\ldots,y_n\}$ such that $(d_F(x_i))=(d_F(y_i))=(6,5,4,3,3,1)$.
However, $(A;B)=(5,4,3,3,3,1;5,5,3,3,3,1)$ is not realizable by a
graph because $(CDZ)$ does not hold for $t=2$.

 {\bf Theorem 2.4.} \ \ Let $A=(a_1,\ldots,a_n)$
and $B=(b_1,\ldots,b_n)$ be two sequences of nonnegative integers in
good order such that $a_i \le b_i$ for $1\le i\le n$. Then the pair
$(A;B)$ is realizable if and only if
$$\sum_{i=1}^t a_i\le t(n-m-1)+\sum_{i=n-m+1}^n
b_i-\varepsilon(t)\mbox{ for each $t$ with $0\le t\le
n$},\eqno{(2)}$$ where $n-m\ge t$.

{\bf Theorem 2.5.} \ \ Let $A=(a_1,\ldots,a_n)$ and
$B=(b_1,\ldots,b_n)$ be two sequences of nonnegative integers in
good order such that $a_i \le b_i$ for $1\le i\le n$. Then the pair
$(A;B)$ is realizable if and only if
$$\sum_{i=1}^t a_i\le \sum_{i=t+1}^n
b_i+\sum_{i=1}^t min\{a_i,t-1\}-\varepsilon(t)\mbox{ for each $t$
with $0\le t\le n.$}\eqno{(3)}$$

{\bf Theorem 2.6.} \ \ Let $A=(a_1,\ldots,a_n)$ and
$B=(b_1,\ldots,b_n)$ be two sequences of nonnegative integers in
good order such that $a_i \le b_i$ for $1\le i\le n$. Then the pair
$(A;B)$ is realizable if and only if
$$\sum_{i=1}^t max\{t-1,a_i\}\le  t(t-1)+\sum_{i=t+1}^n
b_i -\varepsilon(t)\mbox{ for each $t$ with $0\le t\le
n$.}\eqno{(4)}$$

{\bf Theorem 2.7.} \ \ Suppose that $A=(a_1,\ldots,a_n)$ and
$B=(b_1,\ldots,b_n)$ are in good order and  $(b^* _1,\ldots,b^* _n)$
is the conjugate of $(b_1,\ldots, b_n)$. Then the pair $(A;B)$ is
realizable if and only if
$$\sum_{i=1}^t a_i\le \sum_{i=1}^t
(b^*_i-1)-\varepsilon(t)\mbox{ for each $t$ with $0\le t\le
s-1,$}\eqno{(5)}$$ where $s=max\{i|a_i\ge i-1\}$.

\section*{3. Preliminary Results}
\hskip\parindent In this section we present two lemmas, which will
be useful as we proceed with the proofs of Theorem 2.4 and Theorem
2.5.

{\bf Lemma 3.1.} \ \ Let $p_1,\ldots,p_n$ be a sequence of
nonnegative integers. Then for each $t$, $1\le t\le n$,
$$\begin{array}{lll}
  \sum\limits_{i=1}^t p_i+\sum\limits_{i=1}^t max\{-t+1,-p_i\}
  &=&\sum\limits_{i=1}^t max\{0,p_i-t+1\},
 \end{array}$$
 $$\begin{array}{lll}
  \sum\limits_{i=1}^t p_i+\sum\limits_{i=1}^t max\{-t+1,-p_i\}
  &=&\sum\limits_{i=1}^t [max\{t-1,p_i\}-(t-1)].
 \end{array}$$

$\emph Proof.$ We just need to verify the former inequality, as the
latter one follows by similar arguments.

We shall prove the first equality by induction on $t$. The result is
trivially true if $t=1$, since then
$p_1+max\{-p_1,0\}=p_1=max\{p_1,0\}$. Suppose that it holds for all
$t\le k-1$. Let $t=k$, it suffices to show
$p_k+max\{-p_k,-k+1\}=max\{p_k-k+1,0\}$. Suppose that $p_k\ge k-1$,
then $p_k+max\{-p_k,-k+1\}=p_k-k+1=max\{p_k-k+1,0\}$. Otherwise,
$p_k+max\{-p_k,-k+1\}=max\{p_k-k+1,0\}=0$. Thus in each case we have
$p_k+max\{-p_k,-k+1\}=max\{p_k-k+1,0\}$. The result follows by the
principle of induction.

{\bf Lemma 3.2.} \ \ Let $s$ be the largest integer such that
$a_s\ge s-1$. If $A$ and $B$ are in good order, then $(A;B)$ is
realizable if and only if
$$\sum_{i=1}^t a_i\le t(t-1)+\sum_{i=t+1}^n min\{t,b_i\}-\vp(t)\mbox{ for
each $t$ with $0\le t\le s$.}\eqno{(CDZ')}$$

$\emph Proof.$ The lemma states that the number of inequalities to
check in $(CDZ)$ can be reduced, that is, $(CDZ)$ needs not to be
verified for all $t\le n$, but only for each $t\le s$. In the
following, we will show the equivalence between $(CDZ)$ and
$(CDZ')$, which implies that the result of Lemma 3.2 is true. For
$t> s$, we have to consider three cases depending on whether
$b_{t+1}<t,\ldots,b_n<t$ or not.

Case 1. If $b_{t+1}<t,\ldots,b_n<t$, then $(CDZ)$ reduces to
$\sum\limits_{i=1}^t a_i\le t(t-1)+\sum\limits_{i=t+1}^n
b_i-\vp(t).$ Setting $f(t)=t(t-1)+\sum\limits_{i=t+1}^n
b_i-\sum\limits_{i=1}^t a_i-\vp(t)$, so that
$f(t+1)=t(t+1)+\sum\limits_{i=t+2}^n b_i-\sum\limits_{i=1}^{t+1}
a_i-\vp(t+1)$ and hence
$$\begin{array}{lll}
  f(t+1)-f(t)
  &=&2t-b_{t+1}-a_{t+1}+\vp(t)-\vp(t+1)\\
  &\ge& 2(t-b_{t+1})+\vp(t)-\vp(t+1)\\
  &\ge&2(t-b_{t+1})-1\\
  &>&0.
 \end{array}$$

Case 2. If there exists some $i$ with $t<i\le n$ such that $b_i\ge
t$, then rearrange $b_{t+1},\ldots,b_n$ so that they are in
non-increasing order and write $b'_{t+1}\ge\ldots\ge b'_k\ge
t>b'_{k+1}\ge \ldots\ge b'_n$. From this, we obtain
$$\begin{array}{lll}
  f(t)&=&t(t-1)+\sum\limits_{i=t+1}^n
min\{b_i,t\}-\sum\limits_{i=1}^t a_i-\vp(t)\\
  &=&t(t-1)+\sum\limits_{i=t+1}^n
min\{b'_i,t\}-\sum\limits_{i=1}^t a_i-\vp(t)\\
  &=& t(t-1)+\sum\limits_{i=t+1}^k
min\{b'_i,t\}+\sum\limits_{i=k+1}^n
min\{b'_i,t\}-\sum\limits_{i=1}^t a_i-\vp(t)\\
&=&t(t-1)+t(k-t)+\sum\limits_{i=k+1}^n b'_i-\sum\limits_{i=1}^t
a_i-\vp(t)\\
  &=&t(k-1)+\sum\limits_{i=k+1}^n b'_i-\sum\limits_{i=1}^t
  a_i-\vp(t)
 \end{array}$$
and so $f(t+1)=(t+1)(k-1)+\sum\limits_{i=k+1}^n
b'_i-\sum\limits_{i=1}^{t+1} a_i-\vp(t+1)$. Therefore,
$$\begin{array}{lll}
  f(t+1)-f(t)
  &=&k-1-a_{t+1}+\vp(t)-\vp(t+1)\\
  &\ge& k-2-a_{t+1}\\
  &\ge&t-1-a_{t+1}\\
  &>&0,
 \end{array}$$
 the last inequality holds by the fact that the definition of $s$ and $t>s$.

Case 3. If $b_{t+1}\ge t,\ldots,b_n\ge t$, it is a special case of
Case 2, i.e., $k=n$.

Thus in each case we have $f(t+1)>f(t)$. In other words, assuming
$(CDZ')$ holds for all $t=s$, it will also hold for all $s<t\le n$,
since $f(n)\ge \ldots\ge f(s)$. This indicates  the inequality in
Lemma 3.2 holds for all $t\le n$ and thus the equivalence between
Lemma 3.2 and Theorem 1.7. So we are done.

\section*{4. Proofs of Theorem 2.1 and 2.2 }\hskip\parindent
 {\bf Proof of 2.1.}  Since  $(A;B)$ is realizable, there exists a
graph with vertices $v_1,\ldots,v_n$ such that $a_i\le d(v_i)\le
b_i$ for $1\le i\le n$, from which we get $\bar{d}(v_i)\le
\bar{b}_i$ for all $i$ and further that
$\sum\limits_{i=1}^{t}\bar{d}(v_i)\le
\sum\limits_{i=1}^{t}\bar{b}_i$ for each $t\le n$. On the other
hand, by Berge criterion, we know $\sum\limits_{i=1}^{t}d(v_i)\le
\sum\limits_{i=1}^{t}\bar{d}(v_i)$. Therefore
$\sum\limits_{i=1}^{t}a_i\le \sum\limits_{i=1}^{t}d(v_i)$ tells us
that we have $\sum\limits_{i=1}^{t}a_i\le
\sum\limits_{i=1}^{t}\bar{b}_i$ for $1\le t\le n$. \quad $\Box$

{\bf Proof of 2.2.} Suppose $(1')$ holds. Consider the
$(0,1)$-matrix corresponding to $(\bar{b}_1,\ldots,\bar{b}_n)$. Take
any $t\in <n>$. $\sum\limits_{i=1}^t \bar{b}_i$ is the number of 1's
in the first $t$ columns. In this (0,1)-matrix all diagonal elements
are 0, which means that in the sub-matrix consisting of the first
$t$ rows and columns at most $t^2-t$ entries are 1. On the other
hand, each row $j$ has precisely $min\{t,b_j\}$ 1's on the first $t$
points. Hence, for each $t=0,\ldots,n$, we find
$$\sum\limits_{i=1}^t \bar{b}_i\le t(t-1)+\sum\limits_{j=t+1}^n min\{b_j,t\}$$
By the hypothesis, for each $t$ with $0\le t\le n$, we have
$$\sum_{i=1}^t a_i\le t(t-1)+\sum\limits_{j=t+1}^n min\{b_j,t\}
-\varepsilon(t).$$ It then follows from Theorem 1.8 that $(A;B)$ is
realizable. \quad $\Box$

\section*{5. Proof of  Theorem 2.3}\hskip\parindent
Since $(A;B)$ is realizable by a graph, there exists a graph $H$
such that $a_i\le d_H(v_i)\le b_i$ for $1\le i\le n$. Now write
$d_1=d_H(v_1),\ldots,d_n=d_H(v_n)$ and set $f_1=max\{i|d_i\ge i\}$,
$f_2=max\{i|a_i\ge i\}$ and $f=max\{i|b_i\ge i\}$. Note that $f_2\le
f_1\le f$ and
$(\tilde{d_1},\ldots,\tilde{d_n};\tilde{d_1},\ldots,\tilde{d_n})=(d_1+1,\ldots,d_{f_1}+1,d_{f_1+1},
\ldots,d_n;d_1+1,\ldots,d_{f_1}+1,d_{f_1+1}, \ldots,d_n)$ is
bipartite realizable by Ryser criterion. It is then clear that
$\tilde{a}_i\le d_i+1\le \tilde{b}_i$ for $1\le i\le f_1$ and
$\tilde{a}_i\le d_i\le \tilde{b}_i$ for $f_1<i\le n$, which is
equivalent to that $(S;S)$ is interval  realizable. \quad $\Box$

\section*{6. Proofs of Theorem 2.4 to Theorem 2.7}
\hskip\parindent Let $A=(a_1,\ldots,a_n)$ and $B=(b_1,\ldots,b_n)$
be in good order with $a_i\le b_i$ for $1\le i\le n$. Then the
following holds:

 Each of the criteria (2.4)-(2.7) is equivalent to the statement
 that the pair $(A;B)$ is realizable.

$\emph Proof.$ To prove Theorem 2.4 to Theorem 2.7, we go through
the following implications cycle;

$(A;B)$ is realizable $\stackrel{I}\Longrightarrow$ (2)
$\stackrel{II}\Longrightarrow$ (3) $\stackrel{III}\Longrightarrow$
(4) $\stackrel{IV}\Longrightarrow$ (5) $\stackrel{V}\Longrightarrow$
$(CDZ')$$\stackrel{Lem 3.2}\Longrightarrow$ (CDZ)$\stackrel{Thm
1.8}\Longrightarrow$ $(A;B)$ is realizable.

 $\emph I.$   \quad We may assume that $b_i\ge 1$ for $1\le i\le n$.
If $t=0$ and $\varepsilon(0)=1$, then $I_0=<n>$, $a_i=b_i$ for $i\in
<n>$ and $\sum\limits_{i=1}^n b_i\equiv1 (mod2)$. Since $(A;B)$ is
realizable, there exists a graph $G$ such that $a_i\le d_G(v_i)\le
b_i$ for $1\le i\le n$. Thus $\sum\limits_{i=1}^nd(v_i)=
\sum\limits_{i=1}^n b_i\equiv1 (mod2)$, a contradiction, which
proves that $\varepsilon(0)=0$, so that (2) holds in this case. If
$t\neq 0$, then

$$\begin{array}{lll}
  \sum\limits_{i=1}^{t} a_i\le \sum\limits_{i=1}^{t} d_G(v_i)
  &\le&t(t-1)+\sum\limits_{i=t+1}^n min \{d_G(v_i),t\}\\
  &\le&t(t-1)+\sum\limits_{i=t+1}^n min \{b_i,t\}\\
  &=& t(t-1)+\sum\limits_{i=t+1}^{n-m} min \{b_i,t\}+\sum\limits_{i=n-m+1}^n min \{b_i,t\}\\
  &\le& t(t-1)+t(n-m-t)+\sum\limits_{i=n-m+1}^n min \{b_i,t\}\\
  &=& t(n-m-1)+\sum\limits_{i=n-m+1}^n min \{b_i,t\}, \\
 \end{array}$$
 where $n-m\ge t$ and the second inequality because Theorem 1.3.
Thus (1) holds when $\varepsilon(t)=0$. If $\varepsilon(t)=1$, then
$a_i=d_G(v_i)=b_i$ for $i\in I_t$ and $\sum\limits_{i\in
I_t}d_G(v_i)+t|I_t|\equiv 1(mod2)$. Suppose (2) does not hold, then
$$\begin{array}{lll}
\sum\limits_{i=1}^{t} a_i&>&t(n-m-1)+\sum\limits_{i=n-m+1}^n min
  \{b_i,t\}-\varepsilon(t)\\
  &\ge& t(n-m-1)+\sum\limits_{i=n-m+1}^n min \{b_i,t\}-1,\\
\end{array}$$
namely, $$
  \sum\limits_{i=1}^{t} a_i\ge t(n-m-1)+\sum\limits_{i=n-m+1}^n min \{b_i,t\}.
 $$
From the above discussion, we can see that
$$
  \sum\limits_{i=1}^{t} a_i= t(t-1)+\sum\limits_{i=t+1}^n min \{d_G(v_i),t\}.
 $$
Write $J_t=<n>\setminus I_t$, then by the above equality, it follows
that $N(v_i)\subseteq T(i\in J_t)$, each $v_k$ with $k\in I_t$ is
adjacent to each $v_j$ with $j\in <t>$ and not adjacent to any $v_i$
with $i\in J_t$. Consider the subgraph $G_0$ induced by indices in
$I_t$. The sum of the degrees of vertices in $G_0$ is
$\sum\limits_{i\in I_t}d_G(v_i)-t|I_t|=(\sum\limits_{i\in
  I_t}d_G(v_i)+t|I_t|)-2t|I_t|$.
  But this sum is odd since $\sum\limits_{i\in I_t}d_G(v_i)+t|I_t|
  =\sum\limits_{i\in I_t}b_i+t|I_t|\equiv 1(mod 2)$, a contradiction and so the conclusion is valid.

$\emph II.$ \quad Take any $t\le n$. If $a_1,\ldots,a_t< t-1$, then
$\sum\limits_{i=1}^t min \{a_i,t-1\}=\sum\limits_{i=1}^t a_i$ and
thus (3) is trivially true. So we may assume that $a_1,\ldots,a_l\ge
t-1$ and $a_{l+1}, \ldots,a_t<t-1$ for some $l<t$. It follows that
$$\begin{array}{lll}
\sum\limits_{i=1}^t a_i&=&\sum\limits_{i=1}^l a_i+\sum\limits_{i=l+1}^t a_i \quad(using (2))\\
&\le& l(n-m-1)+\sum\limits_{i=n-m+1}^n
b_i-\varepsilon(t)+\sum\limits_{i=l+1}^t a_i \quad(using \quad m=n-t)\\
&=& l(t-1)+\sum\limits_{i=t+1}^n b_i+\sum\limits_{i=l+1}^t
a_i-\varepsilon(t)\\
&=& \sum\limits_{i=1}^l(t-1)+\sum\limits_{i=t+1}^n
b_i+\sum\limits_{i=l+1}^t min\{a_i,t-1\}-\varepsilon(t)\\
&=& \sum\limits_{i=1}^l min\{a_i,t-1\}+\sum\limits_{i=t+1}^n
b_i+\sum\limits_{i=l+1}^t min\{a_i,t-1\}-\varepsilon(t)\\
&=& \sum\limits_{i=t+1}^n b_i+\sum\limits_{i=1}^t min\{a_i,t-1\}-\varepsilon(t),\\
\end{array}$$ as claimed.

$\emph III.$ Assuming (3) holds, an application of Lemma 3.1 gives
 $$\begin{array}{lll}
\qquad\sum\limits_{i=1}^t a_i-\sum\limits_{i=1}^t min\{a_i,t-1\}\le
\sum\limits_{i=t+1}^n b_i-\varepsilon(t) \\
\Longleftrightarrow\sum\limits_{i=1}^t a_i+\sum\limits_{i=1}^t
max\{-a_i,-t+1\}\le
\sum\limits_{i=t+1}^n b_i-\varepsilon(t) \\
\Longleftrightarrow\sum\limits_{i=1}^t max\{0,a_i-t+1\}\le
\sum\limits_{i=t+1}^n b_i-\varepsilon(t) \\
\Longleftrightarrow\sum\limits_{i=1}^t [max\{t-1,a_i\}-(t-1)]\le
\sum\limits_{i=t+1}^n b_i-\varepsilon(t) \\
\Longleftrightarrow\sum\limits_{i=1}^t max\{t-1,a_i\}\le
t(t-1)+\sum\limits_{i=t+1}^n b_i-\varepsilon(t), \\
\end{array}$$ as required.

$\emph IV.$ \quad Suppose to the contrary that there is an index
$t\le s-1$ such that $\sum\limits_{i=1}^t a_i> \sum\limits_{i=1}^t
(b^*_i-1)-\varepsilon(t).$ Then by the definition of $s$, we get
$b_1,\ldots,b_t,\ldots,b_{s-1},b_s\ge t$. Note that the order of the
terms $b_{t+1},\ldots,b_n$ in $B=(b_1,\ldots,b_n)$ clearly has no
bearing on the result of Theorem 2.5, and we shall find it
convenient to assume that $b_{t+1},\ldots,b_n$ are rearranged in
non-increasing order and denote $ b'_{t+1}\ge\ldots \ge b'_n$. For
convenience, we write $b_1=b'_1,\ldots,b_t=b'_t$. Notice that
$b'_i\ge b_i\ge a_i$ for $i\in <n>$ and $A=(a_1,\ldots,a_n)$ and
$B'=(b'_1,\ldots,b'_n)$ are in good order. There now exists an index
$m\in \{t+1,\ldots,n\}$ such that $b'_i\ge t$ for $s\le i\le m$ and
$b'_i<t$ for $m<i\le n$.

Let $B^*$ denote the (0,1)-matrix in which the $b_i$ leading terms
in each row $i$ are ones and the remaining entries are zeros, and
let $b^*_1,\ldots,b^*_n$ be the column sums of $B^*$. Now consider
the left matrix $B^*$ below and the sum $\sum\limits_{i=1}^t b^*_i$.
Row $i(1\le i\le m)$ of $B^*$ contributes $t$ to this sum and row
$i(m<i\le n)$ of $B^*$ contributes $b'_i$ to this sum.  We then have
$$\begin{array}{lll}
\sum\limits_{i=1}^t a_i &>& \sum\limits_{i=1}^t
(b^*_i-1)-\varepsilon(t)\\
&=&\sum\limits_{i=1}^t b^*_i-t-\varepsilon(t)\\
&=&mt+\sum\limits_{i=m+1}^n b'_i-t-\varepsilon(t)\\
&\ge&(m-1)t+\sum\limits_{i=1}^m
max\{a_i,m-1\}-m(m-1)\\
&=&(m-1)(t-m)+\sum\limits_{i=1}^m max\{a_i,m-1\},\end{array}$$
where, in the fourth line, we made use of $(4)$ since
$A=(a_1,\ldots,a_n)$ and $B'=(b'_1,\ldots,b'_n)$ are in good order
and $a_i\le b'_i$ for $i\in <n>$. Hence,
$$\begin{array}{lll}
 \sum\limits_{i=1}^m
max\{a_i,m-1\}&<&(m-1)(m-t)+\sum\limits_{i=1}^t a_i\\
&=&\sum\limits_{i=1}^t a_i+\sum\limits_{i=t+1}^m (m-1)(\mbox{ with
 $\sum\limits_{i=m+1}^m (m-1)$ defined 0})\\
 &\le&\sum\limits_{i=1}^t max\{a_i,m-1\}+\sum\limits_{i=t+1}^m
max\{a_i,m-1\}\\
&=&\sum\limits_{i=1}^m max\{a_i,m-1\}.\end{array}$$ This is a
contradiction, which implies  (5) holds.

$$\left({\begin{array}{c|ccc}
{{1}}\ldots\ldots{{1}}&{{1}}\ldots\ldots{{1}}\\
\ldots&\ldots\\
 {{1}}\ldots\ldots{{1}}&{{1}}\ldots\ldots{{0}}\\\hline
{{1}}\ldots\ldots{{0}}&{{0}}\ldots\ldots{{0}}\\
\ldots&\ldots\\
{{1}}\ldots\ldots{{0}}&{{0}}\ldots\ldots{{0}}\\
\end{array}}\right)\quad\quad
\left({\begin{array}{c|ccc}
{{1}}\ldots\ldots{{1}}&{{1}}\ldots\ldots{{0}}\\
\ldots&\ldots\\
 {{1}}\ldots\ldots{{1}}&{{0}}\ldots\ldots{{0}}\\\hline
{{1}}\ldots\ldots{{0}}&{{0}}\ldots\ldots{{0}}\\
\ldots&\ldots\\
{{1}}\ldots\ldots{{1}}&{{1}}\ldots\ldots{{0}}\\
\end{array}}\right)\quad\quad
\left({\begin{array}{c|ccc}
{{1}}\ldots\ldots{{1}}&{{1}}\ldots\ldots{{0}}\\
\ldots&\ldots\\
 {{1}}\ldots\ldots{{0}}&{{0}}\ldots\ldots{{0}}\\\hline
{{1}}\ldots\ldots{{0}}&{{0}}\ldots\ldots{{0}}\\
\ldots&\ldots\\
{{1}}\ldots\ldots{{1}}&{{1}}\ldots\ldots{{0}}\\
\end{array}}\right)$$

$\emph V.$ \quad For $t<s$, if $b_1\ge\ldots\ge b_t$, one can see
that $b_t\ge a_t\ge a_s\ge s-1\ge t$ and then from the middle matrix
, we deduce
$$\begin{array}{lll}
\sum\limits_{i=1}^t a_i &\le& \sum\limits_{i=1}^t
(b^*_i-1)-\varepsilon(t)\\
&=&t^2+\sum\limits_{i=t+1}^n min\{t,b_i\}-t-\varepsilon(t)\\
&=&t(t-1)+\sum\limits_{i=t+1}^n
min\{t,b_i\}-\varepsilon(t).\end{array}$$ Otherwise, if
$b_1\ge\ldots\ge b_t$ does not hold, then the right matrix gives
$$\begin{array}{lll} \sum\limits_{i=1}^t a_i
&\le&\sum\limits_{i=1}^t
(b^*_i-1)-\varepsilon(t)\\
&<&t^2+\sum\limits_{i=t+1}^n
min\{t,b_i\}-t-\varepsilon(t)\\
&=&t(t-1)+\sum\limits_{i=t+1}^n
min\{t,b_i\}-\varepsilon(t).\end{array}$$ Thus in either case,
$(CDZ')$ is satisfied for all $t<s$. Additionally, the second case
implies that $(CDZ')$ also holds for $t=s$. Hence we conclude that
$(CDZ')$ is true and so the proof is complete. \quad $\Box$


\vskip 0.1cm

\noindent{\bf References}

\vskip 0.1cm

 \bref{[1]}{C. Berge, Graph and Hypergraphs,} {North Holland, Amsterdam, and American
Elsevier,New York}{ (1973).}

\vskip 0.1cm

 \bref{[2]}{B.
Bollob$\acute{a}$s,  Extremal Graph Theory,} { Academic Press, New
York}{ (1978).}

\vskip 0.1cm

 \bref{[3]}{M.C. Cai, X.T. Deng and W.A. Zang, Solution
to a problem on degree sequences of graphs,} { Discrete Math.,}{
219}{ (2000) 253--257.}

\bref{[4]}{M.C. Cai, and L.Y. Kang, A characterization of
box-bounded degree sequences of graphs,} { Graph Comb.,}{ 34}{
(2018) 599--606.}

\vskip 0.1cm

\bref{[5]}{P. Erd\H{o}s and T. Gallai, Graphs with given degrees of
vertices,} { Math. Lapok,}{ 11}{ (1960) 264--274.}

\vskip 0.1cm

 \bref{[6]}{D. R. Fulkerson, A. J. Hofman  and M. H. Mcandrew,
Some properties of graphs with multiple edges,} { Can. J. Math.,}{
17}{ (1965), 166--177.}

\vskip 0.1cm

\bref{[7]}{A. Garg, A. Goel and A. Tripathi,  Constructive
extensions of two results on graphic sequences,} { Discrete Appl.
Math.,}{ 159}{ (2011), 2170--2174.}

\vskip 0.1cm

\bref{[8]}{B. Gr$\ddot{u}$nbaum, Graphs and complexes,} { Report of
the University of Washington, Seattle, Math.,}{ 572B}{ (1969),
(private communication).}

\vskip 0.1cm

\bref{[9]}{J. Y. Guo and J. H. Yin, A variant of Niessen's problem
on degree sequences of graphs,} { Discrete Math. Ther. Comp. Scie,}{
16}{ (2014), 287--292.}

\vskip 0.1cm

\bref{[10]}{W. H$\ddot{a}$sselbarth,  Verzweigtheitvon,} { Match.}{
16}{ (1984), 3--17.}

\vskip 0.1cm

\bref{[11]}{L. Lov\'{a}sz, The factorization of graphs II,} { Acta
Math. Sci. Hungar,}{ 23}{ (1972) 223--246.}

\vskip 0.1cm

\bref{[12]}{T. Niessen, Problem 297 (Research problems),} { Discrete
Math.,}{ 191}{ (1998) 250.}

\vskip 0.1cm

\bref{[13]}{G. Sierksma and H. Hoogeveen, Seven criteria for integer
sequences being graphic, } {J. Graph Theory,} { 15}{ (1991),
223--231.}

\end{document}